\documentclass{article}%
\usepackage{amsmath}
\usepackage{amsfonts,amsthm}
\usepackage{amssymb}
\usepackage{graphicx}
\usepackage{algorithm}
\usepackage{algorithmic}
\usepackage{fancybox}
\usepackage{xcolor}
\usepackage{tikz}
\usepackage{float}%

\newtheorem{theorem}{Theorem}[section]
\newtheorem{corollary}[theorem]{Corollary}
\newtheorem{lemma}[theorem]{Lemma}

\newtheorem{definition}[theorem]{Definition}

\newtheorem{example}[theorem]{Example}
\def\dj{d\kern-0.4em\char"16\kern-0.1em}

\evensidemargin=\oddsidemargin
\begin{document}

\title{The $H$-join of arbitrary families of graphs}
\date{}
\author{Domingos M. Cardoso\thanks{Center for Research and Development in Mathematics and Applications, Department of Mathematics, University of Aveiro, Campus de Santiago, 3810-193 Aveiro, Portugal.
E-mail: dcardoso@ua.pt ORCID: http://orcid.org/0000-0001-6239-3557}
\and Helena Gomes\thanks{Center for Research and Development in Mathematics and Applications, Escola Superior de Educa\c{c}\~{a}o, Instituto Polit\'ecnico de Viseu, Viseu, Portugal. E-mail: hgomes@ua.pt}
\and Sofia J. Pinheiro\thanks{Center for Research and Development in Mathematics and Applications, Department of Mathematics, University of Aveiro, Campus de Santiago, 3810-193 Aveiro, Portugal. E-mail: spinheiro@ua.pt}}

\maketitle

\begin{abstract}
The $H$-join of a family of graphs $\mathcal{G}=\{G_1, \dots, G_p\}$, also called the generalized composition,
$H[G_1, \dots, G_p]$, where all graphs are undirected, simple and finite, is the graph obtained by replacing each
vertex $i$ of $H$ by $G_i$ and adding to the edges of all graphs in $\mathcal{G}$ the edges of the join $G_i \vee G_j$,
for every edge $ij$ of $H$. Some well known graph operations are particular cases of the $H$-join of a family of
graphs $\mathcal{G}$ as it is the case of the lexicographic product (also called composition) of two graphs $H$
and $G$, $H[G]$. During long time the known expressions for the determination of the entire spectrum of the
$H$-join in terms of the spectra of its components and an associated matrix were limited to families of regular
graphs. In this work, we extend such a determination, as well as the determination of the characteristic polynomial,
to families of arbitrary graphs. From the obtained results, the eigenvectors of the adjacency matrix of the $H$-join
can also be determined in terms of the adjacency matrices of the components and an associated matrix.

\medskip


\smallskip


\end{abstract}

\section{Introduction}

Nearly five decades since the publication in 1974 of Allen Shweenk's article \cite{Schwenk1974}, the determination
of the spectrum and characteristic polynomial of the generalized composition $H[G_1, \dots, G_p]$ (recently
designated $H$-join of $\mathcal{G}=\{G_1, \dots, G_p\}$ \cite{Cardoso_et_al2013}), in terms of the spectra
(and characteristic polynomials) of the graphs in $\mathcal{G}$ and an associated matrix, where all graphs
are undirected, simple and finite, was limited to families of regular graphs. Very recently in
\cite{SaravananMuruganArunkumar2020}, as an application of a new generalization of Fiedler's Lemma, the
characteristic polynomial of the universal adjacency matrix of the $H$-joim of a family of arbitrary graphs
is determined in terms of the characteristic polynomial and a related racional function of each component, and
the determinant of an associated matrix. In this work, using a distincy approach, the determination of the spectrum
(and characteristic polynomial) of the $H$-join, in terms of its components and associated matrix, is extended to
families of arbitrary graphs (which should be undirected, simple and finite).\\

The generalized composition $H[G_1, \dots, G_p]$, introduced in \cite[p. 167]{Schwenk1974} was rediscovered in
\cite{Cardoso_et_al2013} under the designation of $H$-join of a family of graphs $\mathcal{G}=\{G_1, \dots, G_p\}$,
where $H$ is a graph of order $p$. In \cite[Th. 7]{Schwenk1974}, assuming that $G_1, \dots, G_p$ are all regular graphs
and taking into account that $V(G_1) \cup \dots \cup V(G_p)$ is an equitable partition $\pi$, the characteristic
polynomial of $H[G_1, \dots, G_p]$ is determined in terms of the characteristic polynomials of the graphs
$G_1, \dots, G_p$ and the matrix associated to $\pi$. Using a generalization of a Fiedler's result
\cite[Lem. 2.2]{Fiedler1974} obtained in \cite[Th. 3]{Cardoso_et_al2013}, the spectrum of the $H$-join of a family of
regular graphs (not necessarily connected) is determined in \cite[Th. 5]{Cardoso_et_al2013}.

When the graphs of the family $\mathcal{G}$ are all isomorphic to a fixed graph $G$, the  $H$-join of $\mathcal{G}$
is the same as the lexicographic product (also called the composition) of the graphs $H$ and $G$ which is denoted as
$H[G]$ (or $H \circ G$). The lexicographic product of two graphs was introduced by Harary in \cite{harary1} and
Sabidussi in \cite{sabidussi} (see also \cite{harary2, hammack_et_al}). From the definition, it is immediate that
this graph operation is associative but not commutative.

In \cite{abreu_et_al}, as an application of the $H$-join spectral properties, the lexicographic powers of a graph
$H$ were considered and their spectra determined, when $H$ is regular. The $k$-th lexicographic power of $H$, $H^k$,
is the lexicographic product of $H$ by itself $k$ times (then $H^2=H[H], H^3=H[H^2]=H^2[H], \dots$). As an example,
in \cite{abreu_et_al}, the spectrum of the $100$-th lexicographic power of the Petersen graph, which has a gogool
number (that is, $10^{100}$) of vertices, was determined. With these powers, $H^k$, in \cite{Cardoso_et_al2017} the
lexicographic polynomials were introduced and their spectra determined, for connected regular graphs $H$, in terms
of the spectrum of $H$ and the coefficients of the polynomial.

Other particular $H$-join graph operations appear in the literature under different designations, as it is the case
of the mixed extension of a graph $H$ studied in \cite{haemers_et_al}, where special attention is given to the mixed
extensions of $P_3$. The mixed extension of a graph $H$, with vertex set $V(H)=\{1, \dots, p\}$, is the $H$-join of
a family of graphs $\mathcal{G}=\{G_1, \dots, G_p\}$, where each graph $G_i \in \mathcal{G}$ is a complete graph or
its complement. From the $H$-join spectral properties, we may conclude that the mixed extensions of a graph $H$ of
order $p$ has at most $p$ eigenvalues unequal to $0$ and $-1$.\\

The remaining part of the paper is organized as follows. The focus of Section~\ref{sec_2} is the preliminaries.
Namely, the notation and basic definitions, the main spectral results of the $H$-join graph operation and the more
relevant properties, in the context of this work, of the main characteristic polynomial and walk-matrix of a graph.
In Section~\ref{sec_3}, the main result of this artice, the determination of the spectrum of the $H$-join
of a family of arbitrary graphs is deduced. Section~\ref{sec_4} includes some final remarks. Namely, about
particular cases of the $H$-join such as the lexicographic product and determination of the eigenvectors of the
adjacency matrix of the $H$-join in terms of the eigenvectors of the adjacency matrices of the components and an
associated matrix.

\section{Preliminaries}\label{sec_2}

\subsection{Notation and basic definitions}

Throughout the text we consider undirected, simple and finite graphs, which are just called graphs. The vertex set
and the edge set of a graph $G$ is denoted by $V(G)$ and $E(G)$, respectively. The order of $G$ is the cardinality
of its vertex set and when it is $n$ we consider that $V(G) = \{1, \dots, n\}$. The eigenvalues of the adjacency
matrix of a graph $G$, $A(G)$, are also called the eigenvalues of $G$.
For each distinct eigenvalue $\mu$ of $G$, ${\mathcal E}_G(\mu)$ denotes the eigenspace of $\mu$ whose dimension is
equal to the algebraic multiplicity of $\mu$, $m(\mu)$. The spectrum of a graph $G$ of order $n$ is denoted by
$\sigma(G)=\{\mu^{[m_1]}_1,\dots,\mu^{[m_s]}_s,\mu^{[m_{s+1}]}_{s+1}, \dots, \mu_t^{[m_t]}\}$, where
$\mu_1,\dots,\mu_s,\mu_{s+1}, \dots, \mu_t$ are the distinct eigenvalues of $G$, $\mu_i^{[m_i]}$ means that
$m(\mu_i)=m_i$ and then $\sum_{j=1}^{t}{m_j}=n$. When we say that $\mu$ is an eigenvalue of $G$ with zero multiplicity
(that is, $m(\mu)=0$) it means that $\mu \not \in \sigma(G)$. The distinct eigenvalues of $G$ are indexed in such way
that the eigenspaces ${\mathcal E}_G(\mu_i)$, for $1 \le i \le s$, are not orthogonal to ${\bf j}_n$, the all-$1$
vector with $n$ entries (sometimes we simple write ${\bf j})$. The eigenvalues $\mu_i$, with $1 \le i \le s$, are
called main eigenvalues of $G$ and the remaining distinct eigenvalues non-main. The concept of main (non-main)
eigenvalue was introduced in \cite{cvetkov70} and further investigated in several publications. As it is well known,
the largest eigenvalue of a connected graph $G$ is main and, when $G$ is regular, all its remaining distinct
eigenvalues are non-main \cite{cds79}. A survey on main eigenvalues was published in \cite{rowmain}.

\subsection{The $H$-join operation}

Now we recall the definition of the $H$-join of a family of graphs \cite{Cardoso_et_al2013}.

\begin{definition}\label{def_h-join}
Consider a graph $H$ with vertex subset $V(H)=\{1, \dots, p\}$ and a family of graphs
$\mathcal{G} = \{G_1, \dots, G_p\}$ such that $|V(G_1)|=n_1, \dots , |V(G_p)|=n_p$.
The $H$-join of $\mathcal{G}$ is the graph
$$
G = \bigvee_{H}{\mathcal{G}}
$$
in which $V(G) = \bigcup_{j=1}^{p}{V(G_j)}$ and
$E(G) = \left(\bigcup_{j=1}^{p}{E(G_j)}\right) \cup \left(\bigcup_{rs \in E(H)}{E(G_r \vee G_s)}\right)$,
where $G_r \vee G_s$ denotes the join.
\end{definition}

\begin{theorem}\label{H-Join_Spectra} \cite{Cardoso_et_al2013}
Let $G$ be the $H$-join as in Definition~\ref{def_h-join}, where $\mathcal{G}$ is a family of regular graphs
such that $G_1$ is $d_1$-regular, $G_2$ is $d_2$-regular, $\dots$ and  $G_p$ is $d_p$-regular. Then
\begin{equation}\label{h-join-spectra}
\sigma(G) = \left(\bigcup_{j=1}^{p}{\left(\sigma(G_j) \setminus \{d_j\}\right)}\right) \cup \sigma(\widetilde{C}),
\end{equation}
where the matrix $\widetilde{C}$ has order $p$ and is such that
\begin{equation}\label{matrix_c}
\left(\widetilde{C}\right)_{rs} = \left\{\begin{array}{ll}
                                          d_r          & \hbox{if } r=s,\\
                                         \sqrt{n_rn_s} & \hbox{if } rs \in E(H),\\
                                          0            & \hbox{otherwise,} \\
                                         \end{array}\right.
\end{equation}
and the set operations in \eqref{h-join-spectra} are done considering possible repetitions of elements of the multisets.
\end{theorem}

From the above theorem, if there is $G_i \in \mathcal{G}$ which is disconnected, with $q$ components, then  its
regularity $d_i$ appears $q$ times in the multiset $\sigma(G_i)$. Therefore, according to \eqref{h-join-spectra},
$d_i$ remains as an eigenvalue of $G$ with multiplicity $q-1$.

From now on, given a graph $H$, we consider the following notation:
$$
\delta_{i,j}(H) = \left\{\begin{array}{ll}
                          1 & \hbox{if } ij \in E(H), \\
                          0 & \hbox{otherwise.}
                          \end{array}
                          \right.
$$

Before the next result, it is worth observe the following. Considering a graph $G$, it is always possible to extend
a basis of the eigensubspace associated to a main eigenvalue $\mu_j$, ${\mathcal E}_G(\mu_j) \cap {\bf j}^{\perp}$,
to one of ${\mathcal E}_G(\mu_j)$ by adding an eigenvector, $\hat{\bf u}_{\mu_j}$, which is orthogonal to
${\mathcal E}_G(\mu_j) \cap {\bf j}^{\perp}$ and uniquely determined, without considering its multiplication by a
nonzero scalar. The eigenvector $\hat{\bf u}_{\mu_j}$ is called the main eigenvector of $\mu_j$. The subspace with
basis $\{\hat{\bf u}_{\mu_1}, \dots, \hat{\bf u}_{\mu_s}\}$ is the main subspace of $G$ and is denoted as $Main(G)$.
Note that for each main eigenvector $\hat{\bf u}_{\mu_j}$ of the basis of $Main(G)$,
$\hat{\bf u}_{\mu_j}^T{\bf j}  \ne 0$.

\begin{lemma}\label{main_and_nom-main_eigenvalues-join}
Let $G$ be the $H$-join as in Definition~\ref{def_h-join} and $\mu_{i,j} \in \sigma(G_i)$. Then
$\mu_{i,j} \in \sigma(G)$ with multiplicity
$$
\left\{\begin{array}{ll}
              m(\mu_{i,j})    & \hbox{whether $\mu_{i,j}$ is a non-main eigenvalue of } G_i, \\
              m(\mu_{i,j})-1 & \hbox{whether $\mu_{i,j}$ is a main eigenvalue of } G_i.
 \end{array}\right.
$$
\end{lemma}

\begin{proof}
Denoting $\delta_{i,j} = \delta_{i,j}(H)$, then $\delta_{i,j} {\bf j}_{n_i}{\bf j}^T_{n_j}$ is an $n_i \times n_j$
matrix whose entries are 1 if $ij \in E(H)$ and $0$ otherwise. Then the adjacency matrix of $G$ has the form
$$
A(G) = \left(\begin{array}{ccccc}
        A(G_1)                            & \delta_{1,2}{\bf j}_{n_1} {\bf j}^T_{n_2} & \cdots & \delta_{1,p-1}{\bf j}_{n_1} {\bf j}^T_{n_{p-1}}&\delta_{1,p}{\bf j}_{n_1} {\bf j}^T_{n_p}\\
       \delta_{2,1}{\bf j}_{n_2} {\bf j}^T_{n_1} & A(G_2)                            & \cdots & \delta_{2,p-1}{\bf j}_{n_2} {\bf j}^T_{n_{p-1}}&\delta_{2,p}{\bf j}_{n_2} {\bf j}^T_{n_p}\\
       \vdots                              &   \vdots                          & \ddots & \vdots                               & \vdots \\
       \delta_{p-1,1}{\bf j}_{n_{p-1}}{\bf j}^T_{n_1} &\delta_{p-1,2}{\bf j}_{n_{p-1}} {\bf j}^T_{n_2} & \cdots & A(G_{p-1})&\delta_{p-1,p}{\bf j}_{n_{p-1}} {\bf j}^T_{n_p}\\
       \delta_{p,1}{\bf j}_{n_p}{\bf j}^T_{n_1}       &\delta_{p,2}{\bf j}_{n_p} {\bf j}^T_{n_2}& \cdots &\delta_{p,p-1}{\bf j}_{n_p} {\bf j}^T_{n_{p-1}} & A(G_p)\\
  \end{array}\right).
$$
Let $\hat{\bf u}_{i,j}$ be an eigenvector of $A(G_i)$ associated to an eigenvalue $\mu_{i,j}$ whose sum of its
components is zero (then, $\mu_{i,j}$  is non-main or it is main with multiplicity greater than one). Then,
\begin{equation}\label{eigenvalue-equation}
A(G) \left(\begin{array}{c}
                 0 \\
                 \vdots \\
                 0 \\
                 \hat{\bf u}_{i,j} \\
                 0 \\
                 \vdots \\
                 0
        \end{array}\right) = \left(\begin{array}{c}
                                    \delta_{1,i}\left({\bf j}^T_{n_i}\hat{\bf u}_{i,j}\right){\bf j}_{n_1}  \\
                                    \vdots \\
                                    \delta_{i-1,i}\left({\bf j}^T_{n_i}\hat{\bf u}_{i,j}\right){\bf j}_{n_{i-1}} \\
                                    A(G_i)\hat{\bf u}_{i,j} \\
                                    \delta_{i+1,i}\left({\bf j}^T_{n_i}\hat{\bf u}_{i,j}\right){\bf j}_{n_{i+1}} \\
                                    \vdots \\
                                    \delta_{p,i}\left({\bf j}^T_{n_i}\hat{\bf u}_{i,j}\right){\bf j}_{n_p}
                             \end{array}\right).
\end{equation}
It should be noted that when $\mu_{i,j}$ is main, there are $m(\mu_{i,j})-1$ linear independent eigenvectors
belonging to ${\mathcal E}_G(\mu_{i,j}) \cap {\bf j}^{\perp}$.
\end{proof}

\subsection{The main characteristic polynomial and the walk-matrix}
If $G$ has $s$ distinct main eigenvalues $\mu_1, \dots, \mu_s$, then the main characteristic polynomial of $G$
is the polynomial of degree $s$ \cite{rowmain}
\begin{eqnarray}
m_G(x) &=& \Pi_{i=1}^{s}{(x-\mu_i)} \nonumber\\
       &=& x^s - c_{0} - c_{1}x - \cdots - c_{s-2} x^{s-2} - c_{s-1} x^{s-1}. \label{mcpG}
\end{eqnarray}
As referred in \cite{rowmain} (see also \cite{CvetkovicPetric84}), if $\mu$ is a main eigenvalue of $G$, so is its
algebraic conjugate $\mu^*$ and then the coefficients of $m_G(x)$ are integers. Furthermore, it is worth to recall
the next result which follows from \cite[Th. 2.5]{Teranishi2001} (see also \cite{rowmain}).

\begin{theorem}\cite[Prop. 2.1]{rowmain}\label{minimal_p}
For every polynomial $f(x) \in \mathbb{Q}[x]$, $f(A(G)){\bf j}={\bf 0}$ if and only if $m_G(x)$ divides $f(x)$.
\end{theorem}

In particular, it is immediate that $m_G(A(G)){\bf j}={\bf 0}$. Therefore,
\begin{equation}\label{main_polynomial}
A^s(G){\bf j} = c_{0}{\bf j} + c_{1}A(G){\bf j} + \cdots + c_{s-2} A^{s-2}(G){\bf j} + c_{s-1} A^{s-1}(G){\bf j}.
\end{equation}

Given a graph $G$ of order $n$, let us consider the $n \times k$ matrix \cite{harscwk2, posu}
$$
{\bf W}_{G;k} = \left({\bf j}, A(G) {\bf j}, A^{2}(G){\bf j}, \ldots, A^{k-1}(G){\bf j} \right ).
$$
The vector space spanned by the columns of ${\bf W}_{G;k}$ is denoted as $ColSp{\bf W}_{G;k}$. The matrix
${\bf W}_{G;k}$ with largest integer $k$ such that the dimension of $ColSp{\bf W}_{G;k}$ is equal to $k$, that is,
such that its columns are linear independent, is referred to be the walk-matrix of $G$ and is just denoted as
${\bf W}_G$. Then, as a consequence of Theorem~\ref{minimal_p} and equality \eqref{main_polynomial}, it follows a
theorem which appears in \cite{hagos1}.

\begin{theorem}\cite[Th. 2.1]{hagos1} \label{hagos}
The rank of $W_G$ is equal to the number of main eigenvalues of the graph $G$.
\end{theorem}

From Theorem~\ref{hagos}, we may conclude that the number of distinct main eigenvalues is
$s = \max \{k: \{{\bf j}, A(G){\bf j}, A^2(G){\bf j}, \ldots, A^{k-1}(G){\bf j}\} \text{ is linearly
independent}\}.$ \\

The equality \eqref{main_polynomial} also implies the next corollary.

\begin{corollary} \label{ap}
The $s$-th column of $A(G){\bf W}_G$ is $ A^{s}(G){\bf j} = {\bf W}_G\left(\begin{array}{c}
                                                                                c_{0} \\
                                                                                \vdots \\
                                                                                c_{s-2} \\
                                                                                c_{s-1} \\
                                                                                \end{array} \right ),$
where $c_j$, for $0 \le j \le s-1$, are the coefficients of the main characteristic polynomial $m_G$,
given in \eqref{mcpG}.
\end{corollary}

This corollary allows the determination of the coefficients of the main characteristic polynomial, $m_G$, by solving
the linear system ${\bf W_G \hat{x}} = A^{s}(G){\bf j}$.\\

From \cite[Th. 2.4]{rowmain} we may conclude the following theorem.

\begin{theorem}\label{colspw}
Let $G$ be a graph with adjacency matrix $A(G)$. Then $ColSp{\bf W_G}$ coincides with $Main(G)$. Moreover $Main(G)$
and the vector space spanned by the vectors orthogonal to $Main(G)$, $\left(Main(G)\right)^{\perp}$, are both
$A(G)$--invariant.
\end{theorem}

From the above definitions, if $G$ is a $r$-regular graph of order $n$, since its largest eigenvalue, $r$,
is the unique main eigenvalue, then $m_G(x) = x - r$ and $W_G = \left( {\bf j}_n \right)$.

\section{The spectrum of the $H$-join of a family of arbitrary graphs}\label{sec_3}

Before the main result of this article, we need to define a special matrix ${\bf \widetilde{W}}$ which will be called
the $H$-join associated matrix.

\begin{definition}\label{main_def}
Let $G$ be the $H$-join as in Definition~\ref{def_h-join} and denote $\delta_{i,j}=\delta_{i,j}(H)$. For each
$G_i \in \mathcal{G}$, consider the main characteristic polynomial \eqref{mcpG},
$m_{G_i}(x)=x^{s_i} - c_{i,0} - c_{i,1}x - \cdots - c_{i,s_i-1} x^{s_i-1}$ and its walk-matrix ${\bf W}_{G_i}$.
The $H$-join associated matrix is the $s \times s$ matrix, with $s= \sum_{i=1}^{p}{s_i}$,
$$
\widetilde{\bf W} = \left(\begin{array}{ccccc}
                       {\bf C}(m_{G_1})    & \delta_{1,2}{\bf M}_{1,2}& \dots & \delta_{1,p-1}{\bf M}_{1,p-1}& \delta_{1,p}{\bf M}_{1,p} \\
                        \delta_{2,1}{\bf M}_{2,1}& {\bf C}(m_{G_2})   & \dots & \delta_{2,p-1}{\bf M}_{2,p-1}& \delta_{2,p}{\bf M}_{2,p} \\
                        \vdots             &    \vdots          &\ddots & \vdots                 & \vdots \\
                        \delta_{p,1}{\bf M}_{p,1}& \delta_{p,2}{\bf M}_{p,2}& \dots & \delta_{p,p-1}M_{p,p-1}& {\bf C}(m_{G_p})
                  \end{array}\right), \text{ where}
$$
${\bf C}(m_{G_i}) = \left(\begin{array}{ccccc}
                           0 &   0   &\dots &   0   & c_{i,0} \\
                           1 &   0   &\dots &   0   & c_{i,1} \\
                           0 &   1   &\dots &   0   & c_{i,2} \\
                           \vdots &\vdots &\ddots&\vdots &\vdots\\
                           0 &   0   &\dots &   1   &c_{i,s_i-1} \\
                          \end{array}\right)$ and ${\bf M}_{i,j}=\left(\begin{array}{c}
                                                                   {\bf j}^T_{n_j}{\bf W}_{G_j} \\
                                                                         0    \; \dots \; 0 \\
                                                                   \vdots \; \ddots\; \vdots\\
                                                                         0    \; \dots \; 0 \\
                                                                   \end{array}\right)$, for $1 \le i,j \le p$.
\end{definition}

Note that the ${\bf C}(m_{G_i})$ is the Frobenius companion matrix of the main characteristic polynomial $m_{G_i}$
and ${\bf M}_{i,j}$ is a $s_i \times s_j$ matrix such that its first row is
${\bf j}_{n_i}^T{\bf W}_{G_i} = (N_0^i, N_1^i, \dots, N_{s_i-1}^i)$, where $N_k^i$ is the number of walks of length $k$ in
$G_i$, for $0 \le k \le s_i-1$ (considering $N^i_0=n_i$) and the entries of the remaining rows are equal to zero.

\begin{theorem}\label{main_theorem}
Let $G$ be the $H$-join as in Definition~\ref{def_h-join}, where $\mathcal{G}$ is a family of arbitrary graphs.
If for each graph $G_i$, with $1 \le i \le p$,
\begin{equation}\label{spectrum_Gi}
\sigma(G_i)=\{\mu_{i,1}^{[m_{i,1}]}, \dots, \mu_{i,s_i}^{[m_{i,s_i}]}, \mu_{i,s_i+1}^{[m_{i,s_i+1}]}, \dots, \mu_{i,t_i}^{[m_{i,t_i}]}\},
\end{equation}
where $m_{i,j}=m(\mu_{i,j})$ and $\mu_{i,1}, \dots, \mu_{i,s_i}$ are the main distinct eigenvalues of $G_i$, then
\begin{eqnarray}
\sigma(G) &=& \bigcup_{i=1}^{p}{\{\mu_{i,1}^{[m_{i,1}-1]}, \dots, \mu_{i,s_i}^{[m_{i,s_i}-1]}\}} \cup
              \bigcup_{i=1}^{p}{\{\mu_{i,s_i+1}^{[m_{i,s_i+1}]}, \dots, \mu_{i,t_i}^{[m_{i,t_i}]}\}} \cup
              \sigma({\bf \widetilde{W}}), \label{G_spectrum}
\end{eqnarray}
where the union of multisets is considered with possible repetitions.
\end{theorem}

\begin{proof}
From Lemma~\ref{main_and_nom-main_eigenvalues-join} it is immediate that
$$
\bigcup_{i=1}^{p}{\{\mu_{i,1}^{[m_{i,1}-1]}, \dots, \mu_{i,s_i}^{[m_{i,s_i}-1]}\}} \cup
              \bigcup_{i=1}^{p}{\{\mu_{i,s_i+1}^{[m_{i,s_i+1}]}, \dots, \mu_{i,t_i}^{[m_{i,t_i}]}\}} \subseteq \sigma(G).
$$
So it just remains to prove that $\sigma(\widetilde{{\bf W}}) \subseteq \sigma(G)$.\\

Let us define the vector
\begin{eqnarray}
\hat{\bf v}      &=& \left(\begin{array}{c}
                            \hat{\bf v}_{1} \\
                            \vdots \\
                            \hat{\bf v}_{p} \\
                     \end{array}\right), \text{ such that } \label{vector_v}\\
\hat{\bf v}_{i}  &=& \sum_{k=0}^{s_i-1}{\alpha_{i,k}A^k(G_i){\bf j}_{n_i}}
                      = {\bf W}_{G_i}\hat{\mathbf{\alpha}}_{i},\label{main_vector_Gi}
\end{eqnarray}
where $\hat{\mathbf{\alpha}}_{i} = \left(\begin{array}{c}
                                    \alpha_{i,0} \\
                                    \alpha_{i,1} \\
                                    \vdots      \\
                                    \alpha_{i,s_i-1} \\
                                    \end{array} \right ), \text{ for } 1 \le i \le p.$
From \eqref{main_vector_Gi}, each $\hat{\bf v}_{i} \in Main(G_i)$ and then all vectors $\hat{\bf v}$ defined in
\eqref{vector_v} are orthogonal to the eigenvectors of $A(G)$ in \eqref{eigenvalue-equation}. Moreover,

\begin{equation}\label{main_subspace}
A(G_i)\hat{\bf v}_{i} = A(G_i){\bf W}_{G_i}\hat{\mathbf{\alpha}}_{i} = \sum_{k=0}^{s_i-1}{\alpha_{i,k}A^{k+1}(G_i){\bf j}_{n_i}}, \text{ for } 1 \le i \le p.
\end{equation}

Therefore,

\begin{eqnarray}
A(G)\hat{\bf v} &=& \left(\begin{array}{cccc}
A(G_1)                                         & \delta_{1,2}{\bf j}_{n_1} {\bf j}^T_{n_2}& \cdots & \delta_{1,p}{\bf j}_{n_1} {\bf j}^T_{n_p}\\
\delta_{2,1}{\bf j}_{n_2} {\bf j}^T_{n_1}      & A(G_2)                                   & \cdots & \delta_{2,p}{\bf j}_{n_2} {\bf j}^T_{n_p}\\
       \vdots                                  &   \vdots                                 & \ddots & \vdots \\
\delta_{p,1}{\bf j}_{n_p}{\bf j}^T_{n_1}       &\delta_{p,2}{\bf j}_{n_p} {\bf j}^T_{n_2} & \cdots & A(G_p)\\
  \end{array}\right) \left(\begin{array}{c}
                     \hat{\bf v}_{1}\\
                     \hat{\bf v}_{2}\\
                     \vdots \\
                     \hat{\bf v}_{p}
  \end{array}\right) \nonumber \\
&=& \left(\begin{array}{c}
                   A(G_1)\hat{\bf v}_{1} + \left(\sum_{k \in [p] \setminus \{1\}}{\delta_{1,k}{\bf j}^T_{n_k}\hat{\bf v}_{k}}\right){\bf j}_{n_1}\\
                   A(G_2)\hat{\bf v}_{2} + \left(\sum_{k \in [p] \setminus \{2\}}{\delta_{2,k}{\bf j}^T_{n_k}\hat{\bf v}_{k}}\right){\bf j}_{n_2}\\
                   \vdots \\
                   A(G_p)\hat{\bf v}_{p} + \left(\sum_{k \in [p] \setminus \{p\}}{\delta_{p,k}{\bf j}^T_{n_k}\hat{\bf v}_{k}}\right){\bf j}_{n_p}
         \end{array}\right) \label{walk_1}\\
&=& \left(\begin{array}{c}
                   A(G_1)\hat{\bf v}_{1} + \left(\sum_{k \in [p] \setminus \{1\}}{\delta_{1,k}{\bf j}^T_{n_k}{\bf W}_{G_k}\hat{\mathbf{\alpha}}_{k}}\right){\bf j}_{n_1}\\
                   A(G_2)\hat{\bf v}_{2} + \left(\sum_{k \in [p] \setminus \{2\}}{\delta_{2,k}{\bf j}^T_{n_k}{\bf W}_{G_k}\hat{\mathbf{\alpha}}_{k}}\right){\bf j}_{n_2}\\
                   \vdots \\
                   A(G_p)\hat{\bf v}_{p} + \left(\sum_{k \in [p] \setminus \{p\}}{\delta_{p,k}{\bf j}^T_{n_k}{\bf W}_{G_k}\hat{\mathbf{\alpha}}_{k}}\right){\bf j}_{n_p}
          \end{array}\right), \label{walk_2}
\end{eqnarray}
where \eqref{walk_2} is obtained applying \eqref{main_vector_Gi} in \eqref{walk_1}. Defining
\begin{equation*}
\beta_{i,0}=\sum_{k \in [p] \setminus \{i\}}{\delta_{i,k}{\bf j}^T_{n_k}{\bf W}_{G_k}\hat{\mathbf{\alpha}}_{k}},
           \text{ for } 1 \le i \le p
\end{equation*}
and taking into account \eqref{main_subspace}, the $i$-th row of \eqref{walk_2} can be written as
{\small \begin{eqnarray}
\hspace{-0.5cm}\beta_{i,0}{\bf j}_{n_i} + A(G_i)\hat{\bf v}_{i}  &=& \left(\underbrace{\sum_{k \in [p]\setminus\{i\}}{\delta_{i,k}{\bf j}^T_{n_k}{\bf W}_{G_k}\hat{\bf \alpha}_{k}}}_{\beta_{i,0}}\right){\bf j}_{n_i} + \sum_{k=0}^{s_i-1}{\alpha_{i,k}A^{k+1}(G_i){\bf j}_{n_i}} \nonumber\\
\hspace{-0.5cm}  &=& \beta_{i,0}{\bf j}_{n_i} + \sum_{k=1}^{s_i-1}{\alpha_{i,k-1}A^k(G_i){\bf j}_{n_i}}
                                      + \alpha_{i,s_i-1}A^{s_i}(G_i){\bf j}_{n_i} \label{ith-row_1}\\
\hspace{-0.5cm}  &=& \beta_{i,0}{\bf j}_{n_i} + \sum_{k=1}^{s_i-1}{\alpha_{i,k-1}A^k(G_i){\bf j}_{n_i}}
                                      + \alpha_{i,s_i-1}{\bf W}_{G_i}\left(\begin{array}{c}
                                                                                c_{i,0} \\
                                                                                c_{i,1} \\
                                                                                \vdots \\
                                                                                c_{i,s_i-1} \\
                                                                            \end{array} \right) \label{ith-row_2}
\end{eqnarray}}
\begin{eqnarray}
\qquad \qquad &=& {\bf W}_{G_i}\left(\begin{array}{c}
                        \beta_{i,0} + \alpha_{i,s_i-1}c_{i,0} \\
                        \alpha_{i,0} + \alpha_{i,s_i-1}c_{i,1}\\
                        \vdots \\
                        \alpha_{i,s_i-2} + \alpha_{i,s_i-1}c_{i,s_i-1} \\
                        \end{array}\right). \label{ith-row_3}
\end{eqnarray}
Observe that \eqref{ith-row_2} is obtained applying Corollary~\ref{ap} to \eqref{ith-row_1}. Taking into account the
definition of $\beta_{i,0}$, \eqref{ith-row_3} can be replaced by the expression

{\tiny
$$
\hspace{-1.5cm}{\bf W}_{G_i}\underbrace{\left(\begin{array}{ccccccccccc}
\overbrace{\delta_{i,1}{\bf j}^T_{n_1}{\bf W}_{G_1}}^{s_1\text{ columns}}&\cdots&\overbrace{\delta_{i,{i-1}}{\bf j}^T_{n_{i-1}}{\bf W}_{G_{i-1}}}^{s_{i-1}\text{ columns}}&      0    &   0  &\cdots&  0   &  c_{i,0} & \overbrace{\delta_{i,{i+1}}{\bf j}^T_{n_{i+1}}{\bf W}_{G_{i+1}}}^{s_{i+1}\text{ columns}} &\cdots&\overbrace{\delta_{i,p}{\bf j}^T_{n_p}{\bf W}_{G_p}}^{s_p\text{ columns}}\\
                  {\bf 0}                                          &\cdots& {\bf 0}&      1    &   0  &\cdots&  0   &  c_{i,1}  &{\bf 0}&\cdots&      {\bf 0}    \\
                  {\bf 0}                                          &\cdots&{\bf 0}&      0    &   1  &\cdots&  0   &  c_{i,2}  &{\bf 0}&\cdots&      {\bf 0}    \\
                  \vdots                                           &\ddots&\vdots&   \vdots  &\vdots&\ddots&\vdots&  \vdots   &\vdots&\ddots&      \vdots     \\
                  {\bf 0}                                          &\cdots&{\bf 0}&      0    &   0  &\cdots&  1   &c_{i,s_i-1} &{\bf 0}&\cdots&      {\bf 0}
\end{array}\right)}_{\widetilde{\bf W}_i}\left(\begin{array}{c}
                        \hat{\bf \alpha}_{1}\\
                        \vdots \\
                        \hat{\bf \alpha}_{i-1}\\
                        \hat{\bf \alpha}_{i}\\
                        \hat{\bf \alpha}_{i+1}\\
                        \vdots \\
                        \hat{\bf \alpha}_{p}
                        \end{array}\right)
$$}
which is equivalent to the expression
$$
{\bf W}_{G_i}\underbrace{\left(\begin{array}{ccccccc}
                   \delta_{i,1}{\bf M}_{i,1} & \dots & \delta_{i,i-1}{\bf M}_{i,i-1} & {\bf C}(m_{G_i}) & \delta_{i,i+1}{\bf M}_{i,i+1} & \dots & \delta_{i,p}{\bf M}_{i,p}\\
                        \end{array}\right)}_{\widetilde{\bf W}_i}\left(\begin{array}{c}
                                                                     \hat{\bf \alpha}_{1}\\
                                                                     \vdots \\
                                                                     \hat{\bf \alpha}_{i-1}\\
                                                                     \hat{\bf \alpha}_{i}\\
                                                                     \hat{\bf \alpha}_{i+1}\\
                                                                     \vdots \\
                                                                     \hat{\bf \alpha}_{p}
                                                                     \end{array}\right).
$$
From the above analysis
\begin{eqnarray*}
A(G)\hat{\bf v} & = & \left(\begin{array}{cccc}
                            {\bf W}_{G_1} &    {\bf 0}    & \cdots &  {\bf 0} \\
                            {\bf 0}   & {\bf W}_{G_2} & \cdots &  {\bf 0} \\
                            \vdots    &    \vdots     & \ddots &  \vdots  \\
                            {\bf 0}   &    {\bf 0}    & \cdots & {\bf W}_{G_p}
                            \end{array}\right) \left(\begin{array}{c}
                                                           {\bf \widetilde{W}}_1\\
                                                           {\bf \widetilde{W}}_2\\
                                                           \vdots\\
                                                           {\bf \widetilde{W}}_p
                                                      \end{array}\right)\left(\begin{array}{c}
                                                                              \hat{\bf \alpha}_{1}\\
                                                                              \hat{\bf \alpha}_{2}\\
                                                                              \vdots \\
                                                                              \hat{\bf \alpha}_{p}
                                                                              \end{array}\right)
\end{eqnarray*}
and, according to \eqref{main_vector_Gi},
\begin{eqnarray*}
\hat{\bf v} &=& \left(\begin{array}{cccc}
                           {\bf W}_{G_1} &    {\bf 0}    & \cdots &  {\bf 0} \\
                           {\bf 0}   & {\bf W}_{G_2} & \cdots &  {\bf 0} \\
                           \vdots    &    \vdots     & \ddots &  \vdots  \\
                           {\bf 0}   &    {\bf 0}    & \cdots & {\bf W}_{G_p}
                           \end{array}\right)\left(\begin{array}{c}
                                                       \hat{\bf \alpha}_{1}\\
                                                       \hat{\bf \alpha}_{2}\\
                                                       \vdots \\
                                                       \hat{\bf \alpha}_{p}
                                                       \end{array}\right). \label{vecto_v}
\end{eqnarray*}
Therefore, $A(G)\hat{\bf v} = \rho \hat{\bf v}$ if and only if
\begin{eqnarray}
\underbrace{\left(\begin{array}{cccc}
       {\bf W}_{G_1} &    {\bf 0}    & \cdots &  {\bf 0} \\
       {\bf 0}       & {\bf W}_{G_2} & \cdots &  {\bf 0} \\
           \vdots    &    \vdots     & \ddots &  \vdots  \\
       {\bf 0}       &    {\bf 0}    & \cdots & {\bf W}_{G_p}
      \end{array}\right)}_{{\bf (*)}}\left(\left(\begin{array}{c}
                               {\bf \widetilde{W}}_1\\
                               {\bf \widetilde{W}}_2\\
                                    \vdots          \\
                               {\bf \widetilde{W}}_p
                              \end{array}\right) - \rho I_s\right)\left(\begin{array}{c}
                                                             \hat{\bf \alpha}_{1}\\
                                                             \hat{\bf \alpha}_{2}\\
                                                             \vdots \\
                                                             \hat{\bf \alpha}_{p}
                                                             \end{array}\right) &=& {\bf 0}. \label{main_equality}
\end{eqnarray}
It is immediate that ${\bf \widetilde{W}} = \left(\begin{array}{c}
                                                  {\bf \widetilde{W}}_1\\
                                                  {\bf \widetilde{W}}_2\\
                                                       \vdots          \\
                                                  {\bf \widetilde{W}}_p
                                                  \end{array}\right)$ and since the columns of each matrix
${\bf W}_{G_i}$ are linear independent, the columns of the matrix $(*)$ are also linear independent. Consequently,
\eqref{main_equality} is equivalent to
\begin{equation}
\left({\bf \widetilde{W}} - \rho I_s\right)\hat{\bf \alpha} = {\bf 0}, \label{w_eigenvetor}
\end{equation}
where $\hat{\bf \alpha} = \left(\begin{array}{c}
                                \hat{\bf \alpha}_{1}\\
                                \vdots \\
                                \hat{\bf \alpha}_{p}
                                \end{array}\right).$
Finally, we may conclude that $(\rho,\hat{\bf v})$ is an eigenpair of $A(G)$ if and only if $(\rho, \hat{\bf \alpha})$,
is an eigenpair of the $H$-join associated matrix ${\bf \widetilde{W}}$.
\end{proof}

Before the next corollary of Theorem~\ref{main_theorem}, it is convenient to introduce the notation
$\phi(G)$ and $\phi({\bf A})$ which, from now on, will be used for the characteristic polynomial of a graph $G$
and a matrix ${\bf A}$, respectively.

\begin{corollary}\label{cor_charact_poly}
Let $G$ be the $H$-join as in Definition~\ref{def_h-join} with associated matrix ${\bf \widetilde{W}}$. Assuming
that $\mathcal{G}=\{G_1, \dots, G_p\}$ is  a family of arbitrary graphs for which $\sigma(G_1), \dots, \sigma(G_p)$
are defined as in \eqref{spectrum_Gi} and $m_{G_1}, \dots, m_{G_p}$ are their main characteristic polynomials, then
$$
\phi(G) = \left(\prod_{i=1}^{p}{\frac{\phi(G_i)}{m_{G_i}}}\right)\phi({\bf \widetilde{W}}).
$$
\end{corollary}

\begin{proof}
Since, for $1 \le i \le p$,
$\sigma(G_i)=\{\mu_{i,1}^{[m_{i,1}]}, \dots, \mu_{i,s_i}^{[m_{i,s_i}]}, \mu_{i,s_i+1}^{[m_{i,s_i+1}]}, \dots, \mu_{i,t_i}^{[m_i,t_i]}\}$,
where the first $s_i$ eigenvalues are main, and the roots of $m_{G_i}$ are all simple main eigenvalues of $G_i$,
it is immediate that $\phi(G_i)=m_{G_i}\phi'(G_i)$, where the roots of the polynomial $\phi'(G_i)$ are the eigenvalues
of $G_i$,
$$
\{\mu_{i,1}^{[m_{i,1}-1]}, \dots, \mu_{i,s_i}^{[m_{i,s_i}-1]}, \mu_{i,s_i+1}^{[m_{i,s_i+1}]}, \dots, \mu_{i,t_i}^{[m_i,t_i]}\}.
$$
Therefore, according to \eqref{G_spectrum}, the roots of $\prod_{i=1}^{p}{\frac{\phi(G_i)}{m_{G_i}}}$ are the
eigenvalues in $\sigma(G) \setminus \sigma({\bf \widetilde{W}})$.
\end{proof}


\begin{example}
Consider the graph $H \cong P_3$, the path with three vertices, and the graphs
$K_{1,3}$, $K_2$ and $P_3$ depicted in the Figure~\ref{figura_1}. Then
$\sigma(K_{1,3})=\{\sqrt{3},-\sqrt{3},0^{[2]}\},$ $\sigma(K_2)=\{1,-1\}$,
$\sigma(P_3)=\{\sqrt{2},-\sqrt{2},0\}$
and their main characteristic polynomials are $m_{K_{1,3}}(x) = x^2 - 3$, $m_{K_2}(x) = x - 1$
and $m_{P_3}(x) = x^2 - 2$, respectively.

\begin{figure}[h]
\begin{center}
\unitlength=0.25 mm
\begin{picture}(400,120)(60,60)
%
\put(50,110){\line(2,1){50}}  
\put(50,110){\line(2,-1){50}} 
\put(50,110){\line(1,0){50}}  
\put(150,85){\line(0,1){50}}  
\put(225,110){\line(0,1){25}} 
\put(225,110){\line(0,-1){25}}
%
%
\put(50,110){\circle*{5.7}}  
\put(100,135){\circle*{5.7}} 
\put(100,110){\circle*{5.7}} 
\put(100,85){\circle*{5.7}}  
%
\put(150,135){\circle*{5.7}} 
\put(150,85){\circle*{5.7}}  
%
\put(225,135){\circle*{5.7}} 
\put(225,110){\circle*{5.7}} 
\put(225,85){\circle*{5.7}} 
%
\put(40,110){\makebox(0,0){\footnotesize 1}}
\put(100,145){\makebox(0,0){\footnotesize 2}}
\put(110,110){\makebox(0,0){\footnotesize 3}}
\put(100,75){\makebox(0,0){\footnotesize 4}}
\put(75,50){\makebox(0,0){\footnotesize $K_{1,3}$}}
%
\put(150,145){\makebox(0,0){\footnotesize 5}}
\put(150,75){\makebox(0,0){\footnotesize 6}}
\put(150,50){\makebox(0,0){\footnotesize $K_2$}}
%
\put(225,145){\makebox(0,0){\footnotesize 7}}
\put(235,110){\makebox(0,0){\footnotesize 8}}
\put(225,75){\makebox(0,0){\footnotesize 9}}
\put(225,50){\makebox(0,0){\footnotesize $P_3$}}
%
%
%
%
\put(300,110){\line(2,1){50}}  
\put(300,110){\line(2,-1){50}} 
\put(300,110){\line(1,0){50}}  
\put(400,85){\line(0,1){50}} 
\put(300,110){\line(4,1){100}} 
\put(300,110){\line(4,-1){100}}
\put(350,135){\line(1,0){125}} 
\put(350,135){\line(1,-1){50}} 
\put(350,110){\line(2,1){50}} 
\put(350,110){\line(2,-1){50}} 
\put(350,85){\line(1,0){125}}  
\put(350,85){\line(1,1){50}}   
\put(400,135){\line(3,-1){75}} 
\put(400,135){\line(3,-2){75}} 
\put(400,85){\line(3,1){75}}   
\put(400,85){\line(3,2){75}}   
\put(475,135){\line(0,-1){50}} 
%
\put(300,110){\circle*{5.7}} 
\put(350,135){\circle*{5.7}} 
\put(350,110){\circle*{5.7}}  
\put(350,85){\circle*{5.7}}  
%
\put(400,135){\circle*{5.7}}
\put(400,85){\circle*{5.7}} 
%
\put(475,135){\circle*{5.7}} 
\put(475,110){\circle*{5.7}} 
\put(475,85){\circle*{5.7}}  
%
\put(295,110){\makebox(0,0){\footnotesize 1}}
\put(350,145){\makebox(0,0){\footnotesize 2}}
\put(360,110){\makebox(0,0){\footnotesize 3}}
\put(350,75){\makebox(0,0){\footnotesize 4}}
%
\put(400,145){\makebox(0,0){\footnotesize 5}}
\put(400,75){\makebox(0,0){\footnotesize 6}}
%
\put(475,145){\makebox(0,0){\footnotesize 7}}
\put(485,110){\makebox(0,0){\footnotesize 8}}
\put(475,75){\makebox(0,0){\footnotesize 9}}
\put(400,50){\makebox(0,0){{\footnotesize $G=\bigvee_{P_3}{\{K_{1,3}, K_2, P_3\}}$}}}
\end{picture}
\end{center}
\caption{The $P_3$-join of the family of graphs $K_{1,3}$, $K_2$ and $P_3$.}\label{figura_1}
\end{figure}
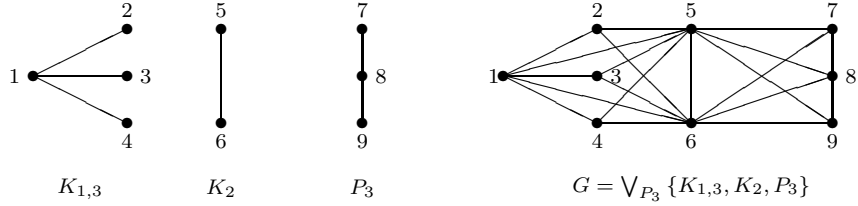
Since
$$
\begin{array}{lclcl}
{\bf \widetilde{W}}_1 &=& \left(\begin{array}{ccccc}
                            0  & c_{1,0} &  \delta_{1,2}2 & \delta_{1,3}3 & \delta_{1,3}4 \\
                            1 & c_{1,1} &         0       &         0     &          0    \\
                           \end{array}\right) & = & \left(\begin{array}{ccccc}
                                                           0  & 3 & 2 & 0 & 0 \\
                                                           1  & 0 & 0 & 0 & 0 \\
                                                          \end{array}\right),\\
{\bf \widetilde{W}}_2 &=& \left(\begin{array}{ccccc}
                           \delta_{2,1}4 & \delta_{2,1}6 & c_{2,0} & \delta_{2,3}3 & \delta_{2,3}4\\
                         \end{array}\right) & = & \left(\begin{array}{ccccc}
                                                         4 & 6 & 1 & 3 & 4\\
                                                        \end{array}\right),\\
{\bf \widetilde{W}}_3 &=& \left(\begin{array}{ccccc}
                           \delta_{3,1}4 & \delta_{3,1}6 & \delta_{3,2}2 &  0  & c_{3,0} \\
                                   0     &          0    &          0    &  1  & c_{3,1}
                         \end{array}\right) & = & \left(\begin{array}{ccccc}
                                                         0 & 0 &  2 &  0 & 2 \\
                                                         0 & 0 &  0 &  1 & 0
                                                        \end{array}\right),
\end{array}
$$
it follows that
\begin{equation*}
{\bf \widetilde{W}} = \left(\begin{array}{c}
                              {\bf \widetilde{W}}_1\\
                              {\bf \widetilde{W}}_2\\
                              {\bf \widetilde{W}}_3
                              \end{array}\right) = \left(\begin{array}{ccccc}
                                                                0 & 3 & 2 & 0 & 0 \\
                                                                1 & 0 & 0 & 0 & 0 \\
                                                                4 & 6 & 1 & 3 & 4 \\
                                                                0 & 0 & 2 & 0 & 2 \\
                                                                0 & 0 & 0 & 1 & 0
                                                          \end{array}\right).
\end{equation*}
Therefore, the characteristic polynomial of ${\bf \widetilde{W}}$ is the polynomial
$$
\phi({\bf \widetilde{W}}) = -42 - 40 x + 15 x^2 + 19 x^3 + x^4 - x^5
$$
and, applying Corollary~\ref{cor_charact_poly}, we obtain the characteristic polynomial of $G$,
$$
\phi(G) = x^3(x+1)\phi({\bf \widetilde{W}}) = x^3(x+1)(-42 - 40 x + 15 x^2 + 19 x^3 + x^4 - x^5).
$$
\end{example}

\section{Final remarks}\label{sec_4}

When all graphs of the family $\mathcal{G}$ are regular, that is, $G_1$ is $d_1$-regular, $G_2$ is $d_2$-regular,
$\dots$, $G_p$ is $d_p$-regular, the walk matrices are ${\bf W}_{G_1}=\left({\bf j}_{n_1}\right)$,
${\bf W}_{G_2}=\left({\bf j}_{n_2}\right)$, $\dots$, ${\bf W}_{G_p}=\left({\bf j}_{n_p}\right)$, respectively.
Consequently, the main polynomials are $m_{G_1}(x) = x - d_1$, $m_{G_2}(x) = x - d_2$, $\dots$, $m_{G_p}(x) = x - d_p$.
As direct
consequence, for this particular case, the $H$-join associated matrix is
{\small $$
{\bf\widetilde{W}}=\left(\begin{array}{cccc}
                          d_1 &  \delta_{1,2}{\bf j}_{n_2}^T{\bf W}_{G_2}  & \cdots & \delta_{1,p}{\bf j}_{n_p}^T{\bf W}_{G_p}\\
                          \delta_{2,1}{\bf j}_{n_1}^T{\bf W}_{G_1} &  d_2  & \cdots & \delta_{2,p}{\bf j}_{n_p}^T{\bf W}_{G_p}\\
                          \vdots                             &\vdots & \ddots & \vdots \\
                          \delta_{p,1}{\bf j}_{n_1}^T{\bf W}_{G_1} & \delta_{p,2}{\bf j}_{n_2}^T{\bf W}_{G_2}  & \cdots &d_p \\
                        \end{array}\right) = \left(\begin{array}{cccc}
                                                    d_1 &  \delta_{1,2}n_2  & \cdots & \delta_{1,p}n_p \\
                                                    \delta_{2,1}n_1 &  d_2  & \cdots & \delta_{2,p}n_p \\
                                                    \vdots          &\vdots & \ddots & \vdots \\
                                                    \delta_{p,1}n_1 &\delta_{p,2}n_2 & \cdots & d_p \\
                                               \end{array}\right).
$$}
Therefore, it is immediate that when all the graphs of the family $\mathcal{G}$ are regular, the matrix
${\bf\widetilde{W}}$ and the matrix $\widetilde{C}$ in \eqref{matrix_c} are similar matrices. Note that
$\widetilde{C} = D {\bf\widetilde{W}}D^{-1}$, where
$D = \text{diag}\left(\sqrt{n_1}, \sqrt{n_2}, \dots, \sqrt{n_p}\right)$ and thus
${\bf\widetilde{W}}$ and $\widetilde{C}$ are cospectral matrices as it should be.\\

In the particular case of the lexicographic product $H[G]$, which is the $H$-join of a family of graphs
$\mathcal{G}$, where all the graphs in $\mathcal{G}$ are isomorphic to a fixed graph $G$, consider that the graph $H$
has order $p$ and the graph $G$ has order $n$. Let
$\sigma(G)=\{\mu_1^{[m_1]}, \dots, \mu_{s}^{[m_s]}, \mu_{s+1}^{[m_{s+1}]}, \dots, \mu_{t}^{[m_t]}\}$, where
$\mu_1, \dots, \mu_s$ are the distinct main eigenvalues of $G$ and $\sum_{i=1}^{t}{m_{i}}=n$. Then, according to
the Definition~\ref{main_def}, the $H$-join associated matrix is
$$
{\bf \widetilde{W}} = \left(\begin{array}{ccccc}
                       {\bf C}(m_{G}) & \delta_{1,2}{\bf M} & \dots & \delta_{1,p-1}{\bf M} & \delta_{1,p}{\bf M} \\
                        \delta_{2,1}{\bf M} & {\bf C}(m_{G})& \dots & \delta_{2,p-1}{\bf M} & \delta_{2,p}{\bf M} \\
                        \vdots        &    \vdots     &\ddots & \vdots          & \vdots \\
                        \delta_{p,1}{\bf M} & \delta_{p,2}{\bf M} & \dots & \delta_{p,p-1}{\bf M} & {\bf C}(m_{G})
                      \end{array}\right),
$$
where ${\bf C}(m_{G})$ is the Frobenius companion matrix of $m_G$ and
${\bf M} = \left(\begin{array}{c}
            {\bf j}_n^{\top}{\bf W}_G\\
            {\bf 0}\\
           \end{array}\right)$ (both are $s \times s$ matrices). Applying Theorem~\ref{main_theorem}, it follows
\begin{eqnarray*}
\sigma(G) &=& p{\{\mu_{1}^{[m_{1}-1]}, \dots, \mu_{s}^{[m_{s}-1]}\}} \cup
              p{\{\mu_{s+1}^{[m_{s+1}]}, \dots, \mu_{t}^{[m_{t}]}\}} \cup
              \sigma({\bf \widetilde{W}}),
\end{eqnarray*}
where the multiplication of $p$ by a set $X$ means the union of $X$ with himself $p$ times. Therefore, from
Corollary~\ref{cor_charact_poly}, the characteristic polynomial of $H[G]$ is
\begin{equation}
\phi(H[G]) = \left(\frac{\phi(G)}{m_G}\right)^p \phi({\bf \widetilde{W}}). \label{charact_poly_lex}
\end{equation}

In \cite[Th. 2.4]{WangWong2018} a distinct expression for $\phi(H[G])$ is determined. Such expression is not
only related with $H$ and $G$ but also with eigenspaces of the adjacency matrix of  $G$.\\

From the obtained results, we are able to determine all the eigenvectors of the adjacency matrix of the $H$-join of
a family of arbitrary graphs $G_1, \dots, G_p$ in terms of the eigenvectors of the adjacency matrices $A(G_i)$, for
$1 \le i \le p$, and the eigenvectors of the $H$-join associated matrix ${\bf \widetilde{W}}$, as follows.
\begin{enumerate}
\item Let $G$ be the $H$-join as in Definition~\ref{def_h-join}, where $\mathcal{G}=\{G_1, G_2, \dots, G_p\}$ is a
      family of arbitrary graphs.
\item For $1 \le i \le p$, consider $\sigma(G_i)$ as defined in \eqref{spectrum_Gi}. For each eigenvalue
      $\mu_{i,j} \in \sigma(G_i)$, every eigenvetor
      $\hat{\bf u}_{i,j} \in \mathcal{E}_{G_i}(\mu_{i,j}) \cap {\bf j}_{n_i}^{\perp}$ defines an eigenvector for
      $A(G)$ as in \eqref{eigenvalue-equation}.
\item The remaining eigenvectors of $A(G)$ are the vectors $\hat{\bf v}$ defined in \eqref{vector_v}-\eqref{main_vector_Gi}
      from the eigenvectors of ${\bf \widetilde{W}}$, $\hat{\bf \alpha}$, obtained as linear independent solutions
      of \eqref{w_eigenvetor}, for each $\rho \in \sigma({\bf \widetilde{W}})$.
\end{enumerate}

\medskip\textbf{Acknowledgments.}
This work is supported by the Center for Research and Development in Mathematics and Applications (CIDMA) through
the Portuguese Foundation for Science and Technology (FCT - Fundação para a Ciência e a Tecnologia), reference
UIDB/04106/2020.

\end{document}